\documentclass[10pt]{article}
\usepackage{amsmath,amssymb,amsthm}
\usepackage{epsfig}
\usepackage{color}

\title{Pancyclicity}

\title{$(r)$-Pancyclic, $(r)$-Bipancyclic and Oddly $(r)$-Bipancyclic Graphs\thanks{Research
supported by NSF REU Grant DMS1262838, University of West Georgia}}
\author{
Abdollah Khodkar\\
Department of Mathematics\\
University of West Georgia\\
Carrollton, GA 30082\\
{\tt akhodkar@westga.edu}\vspace{2mm}\\
Oliver Sawin\\
Department of Mathematics, \\
Rensselaer Polytechnic Institute\\
Troy, NY 12180\\
{\tt osawin@gmail.com}\vspace{2mm}\\
Lisa Mueller\\
Department of Mathematics, \\
California State University, Fullerton\\
Fullerton, CA 92833\\
{\tt exceedinglyhappy@csu.fullerton.edu}\vspace{2mm}\\
WonHyuk Choi\\
Department of Mathematics, Pomona College\\
Claremont, CA 91711\\
{\tt wonhyuk.choi@pomona.edu}\\
}

\date{}

 \newtheorem{theorem}{Theorem}[section]

\newtheorem{conjecture}[theorem]{Conjecture}

\theoremstyle{definition}

\begin{document}

\maketitle
\begin{abstract}
\noindent
A graph with $v$ vertices is $(r)$-pancyclic if
it contains precisely $r$ cycles of every length from 3 to $v$. A bipartite graph with
even number of vertices $v$
is said to be $(r)$-bipancyclic
if it contains precisely $r$ cycles of each even length from 4 to $v$.
A bipartite graph with odd number of vertices $v$ and minimum degree at least 2
is said to be oddly $(r)$-bipancyclic  if it contains precisely $r$ cycles of each even length from 4 to $v-1$.
In this paper, using computer search, we classify all $(r)$-pancyclic and $(r)$-bipancyclic
graphs with $v$ vertices and at most $v+5$ edges. We also classify all
oddly $(r)$-bipancyclic graphs with $v$ vertices and at most $v+4$ edges.
\vspace{2mm}
\end{abstract}

\section{Definitions}
All graphs in this paper are finite, simple and undirected.
A {\em pancyclic} graph, see \cite {Bondy}, of order $v$ is a graph that contains cycles of every
length from 3 to $v$.
Pancyclic graphs
are a generalization of Hamiltonian graphs, which have a
cycle containing all the vertices of the graph.
A pancyclic graph with exactly one cycle of
every possible length is called {\em uniquely pancyclic} (UPC) (see \cite{Shi,Markstr}).

A bipartite graph on $v$ vertices, $v$ even, is called a {\em uniquely bipancyclic} (UBPC) if
it contains precisely one cycle of every even length from 4 to $v$.

The definitions above can be generalized as follows.
A graph with $v$ vertices is {\em $(r)$-pancyclic} \cite{Zam} if
it contains precisely $r$ cycles of every length from 3 to $v$.
Similarly, a bipartite graph with even number of vertices $v$
is said to be {\em $(r)$-bipancyclic}
if it contains precisely $r$ cycles of each even length from 4 to $v$.

A bipartite graph with odd number of vertices $v$ cannot have a cycle of length $v$ and
the largest possible cycle length in such graphs is $v-1$. On the other hand, adding a pendant edge to a graph
does not change the cycle structure of that graph. This bring us to the following definition.
A bipartite graph with odd number of vertices $v$ and minimum degree at least 2
is said to be {\em oddly $(r)$-bipancyclic}
if it contains precisely $r$ cycles of each even length from 4 to $v-1$.

Every $(r)$-pancyclic or $(r)$-bipancyclic graph is Hamiltonian, so we shall represent a $v$-vertex
$(r)$-pancyclic or $(r)$-bipancyclic graph as a Hamilton cycle
together with some edges, which we shall call {\em chords}.
Similarly, a $v$-vertex oddly $(r)$-bipancyclic graph is represented
as a cycle $C$ of length $v-1$ together with some edges, which are either chords or are edges
at the vertex $u\not\in V(C)$, which we shall call {\em chordettes}.

In this paper, using computer search, we classify all $(r)$-pancyclic and $(r)$-bipancyclic
graphs with $v$ vertices and at most $v+5$ edges. We also classify all
oddly $(r)$-bipancyclic graphs with $v$ vertices and at most $v+4$ edges.

\section{Computer search}

The new result presented in the following sections obtained using different computer programs.
The codes for these programs can be found at https://github.com/osawin/Pancyclic.
Each of these programs accepts as input a text file containing a positive integer
to determine the number of cycles of each length and a list of graph schema.
Each schema has a name, and a list of arcs and chords, both represented by a pair of integers denoting their adjacent vertices.
The program goes through the schema one by one looking for valid graphs.
First, the program determinds which edges are arcs, and which are chords or chordettes, based on the formatting of the input.
Then, the program finds all cycles within the schema. Next, the program assigns each arc a variable, representing the size of
the path it will be replaced with to generate a graph. The program tests all possible combinations of values of these variables
to see if it will generate a valid graph. In order to search through graphs faster,
the program tests some classes of combinations of arc lengths. For example, the code assigns all possible combinations of even or
odd to arcs, and examines which combinations will result in a valid
number of even and odd cycle lengths. This allows the code to quickly pass over many non-valid graphs.
Once the program has found all valid graphs, it outputs a description of the schema it was checking,
whether that layout contains any valid graphs,
and the arc lengths that generates valid graphs, if any is found.

\section{$(r)$-pancyclic graphs}

\begin{figure}[ht]
\begin{center}
\begin{tabular}{c}
\includegraphics[scale=.400]{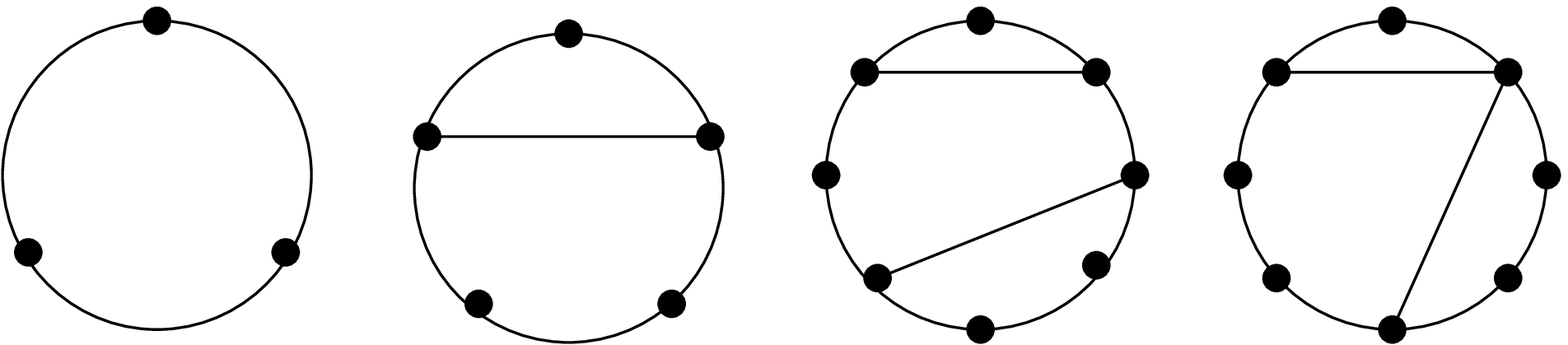}\\
\includegraphics[scale=.400]{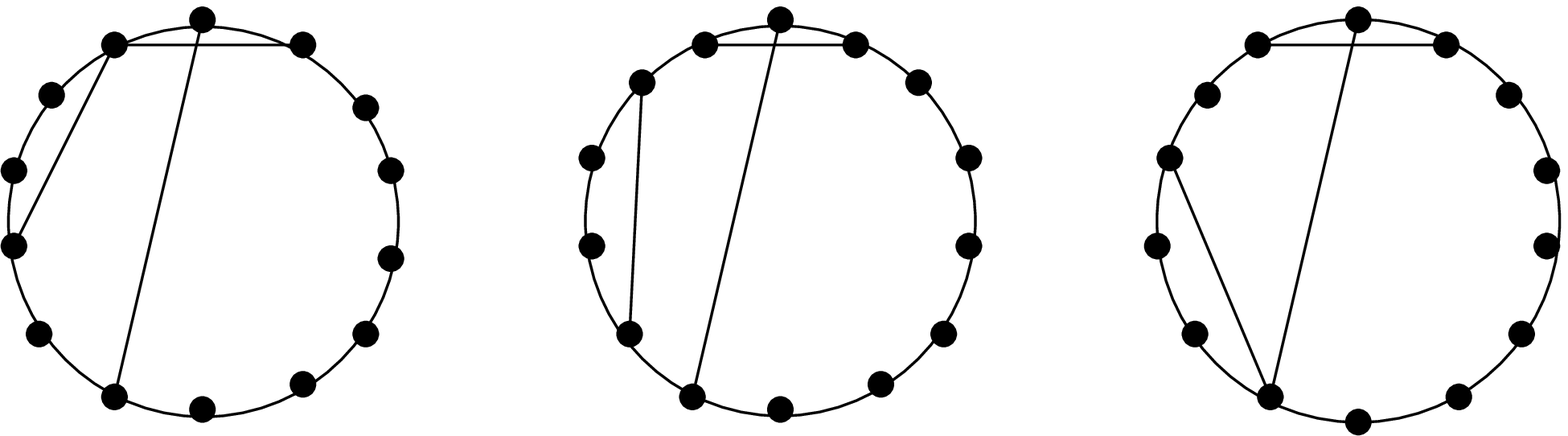}
\end{tabular}
\caption{The (1)-pancyclic graphs with three or fewer chords}
\label{UPC0123chords}
\end{center}
\end{figure}

In \cite{Shi}, Shi classifies all the the UPC graphs with $v$ vertices
and $v+m$ edges for $m=0, 1, 2, 3$ (see Figure \ref{UPC0123chords}) and conjectures that:

\begin{conjecture}\label{Shi.Conj}
There is no UPC graph with $v$ vertices and $v + m$ edges
for $m\geq 4$.
\end{conjecture}

\noindent In \cite{Markstr}, using computer programs, Markstr\"{o}m reconfirms the results
obtained by Shi in \cite{Shi} and shows that there is no UPC graph with $v$ vertices and $v + 5$ edges.

\begin{figure}[ht]
\begin{center}
\includegraphics[scale=.450]{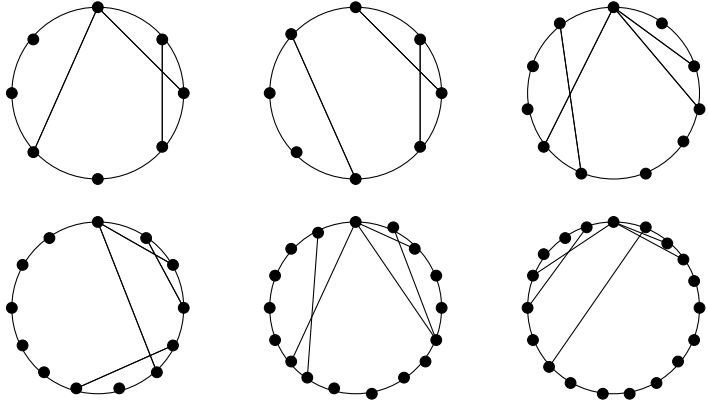}
\caption{The (2)-pancyclic graphs given in \cite{Zam}, $v=8, 11, 13, 17,19$ }
\label{2pan}
\end{center}
\end{figure}

In \cite{Zam}, Zamfirescu studies $(2)$-pancyclic graphs and finds
the two (2)-pancyclic graphs of order 8, which are the smallest (2)-pancyclic graphs.
In the same paper, he proves that there exist no (2)-pancyclic graphs on 9 or 10 vertices.
For each $v\in\{11, 13, 17, 19\}$ a (2)-pancyclic graph of order $v$ is also presented in \cite{Zam}.
(See Figure \ref{2pan}.)

\begin{figure}[ht]
\begin{center}
\begin{tabular}{ccc}
\includegraphics[scale=.380]{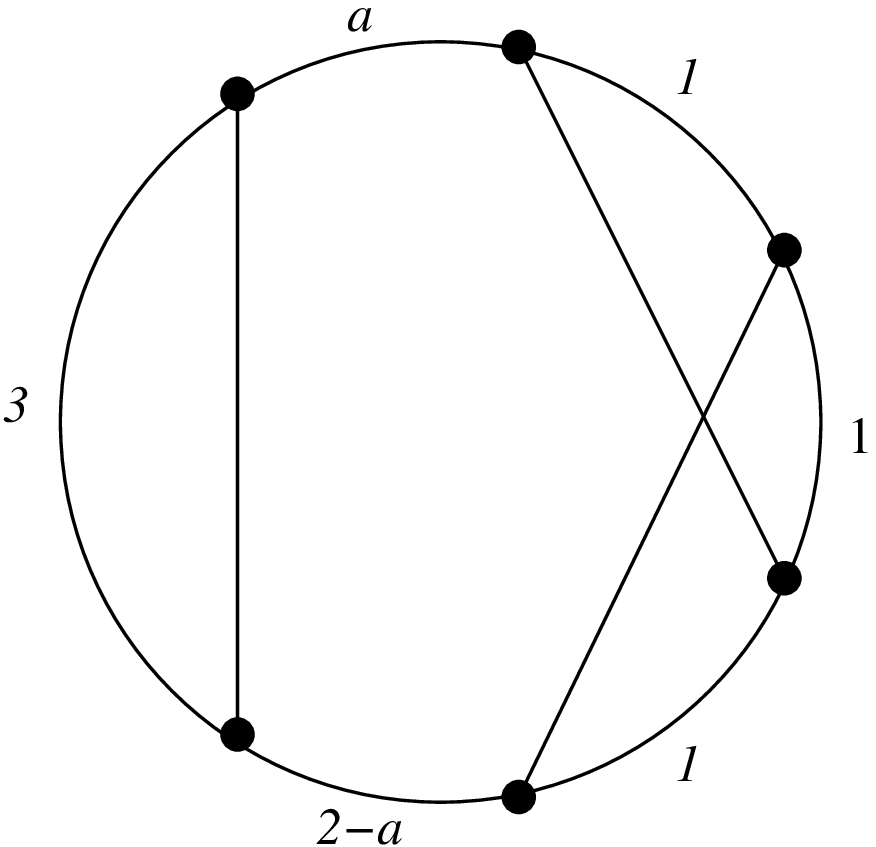}&
\includegraphics[scale=.380]{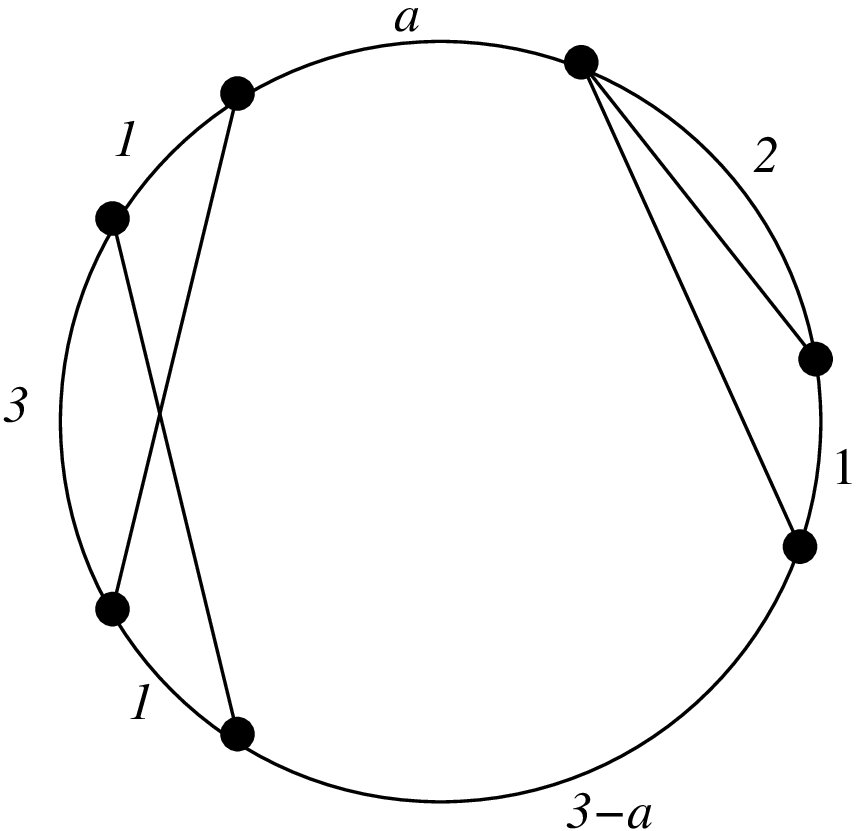}&
\includegraphics[scale=.380]{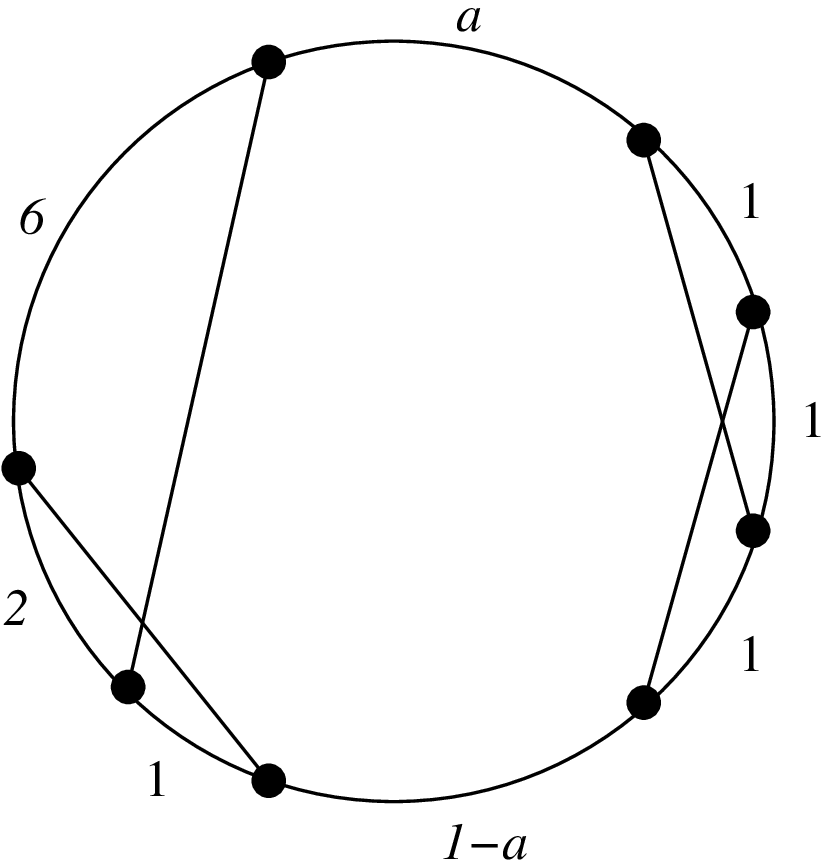}\\
order=8& order=11 & order=13\\
2 non-isomorphic & 4 non-isomorphic & 2 non-isomorphic\\
\includegraphics[scale=.380]{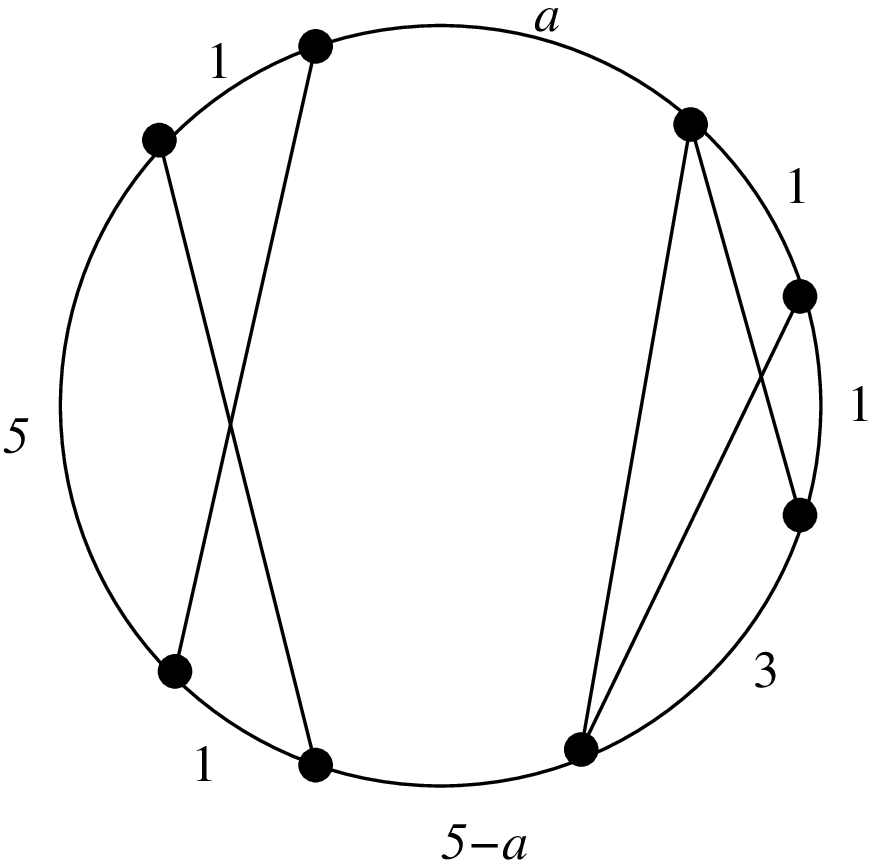}&
\includegraphics[scale=.380]{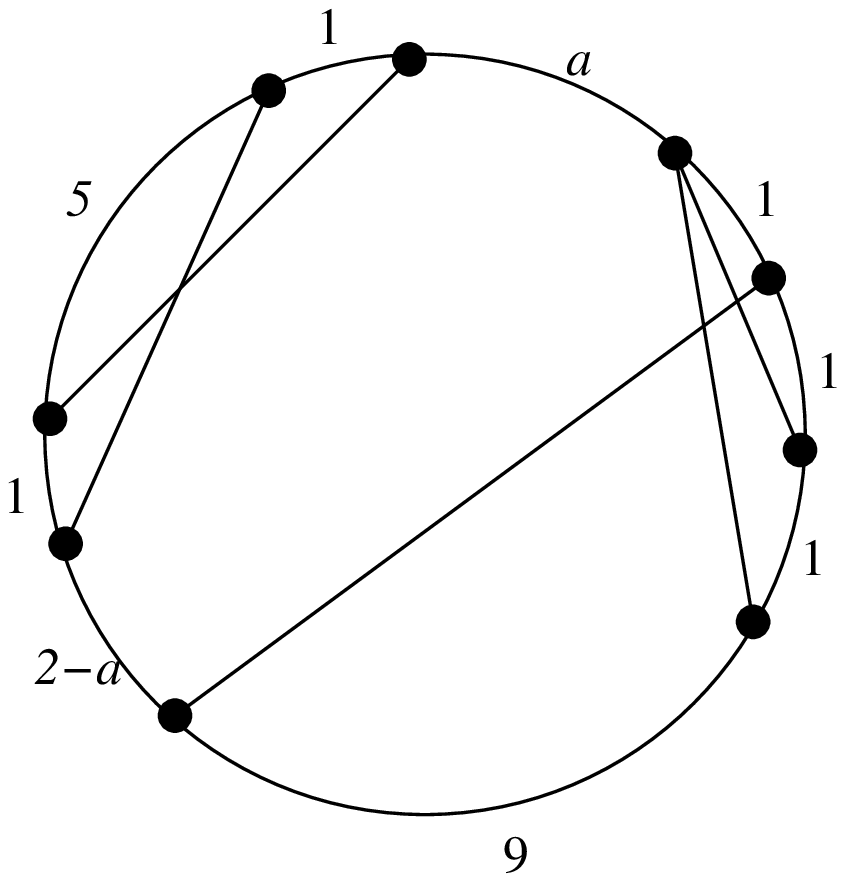}&\\
order=17 & order=19 &\\
6 non-isomorphic & 3 non-isomorphic &\\
\end{tabular}
\caption{The (2)-pancyclic graphs with five chords or less.
The label for an arc indicates the number of edges on that arc.}
\label{all2pancyclic}
\end{center}
\end{figure}

Computer search shows that:

\begin{theorem}\label{(r)-pancyclic}
Let $G$ be a graph with $v$ vertices and $v+m$ edges, where $0\leq m\leq 5$.
\begin{enumerate}
\item The graphs displayed in Figure \ref{all2pancyclic} are the only (2)-pancyclic graphs
with $v$ vertices and $v+m$ edges, where $0\leq m\leq 5$.

\item There is no $(r)$-pancyclic graph, $r\geq 3$, with $v$ vertices and $v+m$ edges,
where $0\leq m\leq 5$.
\end{enumerate}
\end{theorem}

\section{$(r)$-bipancyclic graphs}
\begin{figure}[ht]
\begin{center}
\includegraphics[scale=.350]{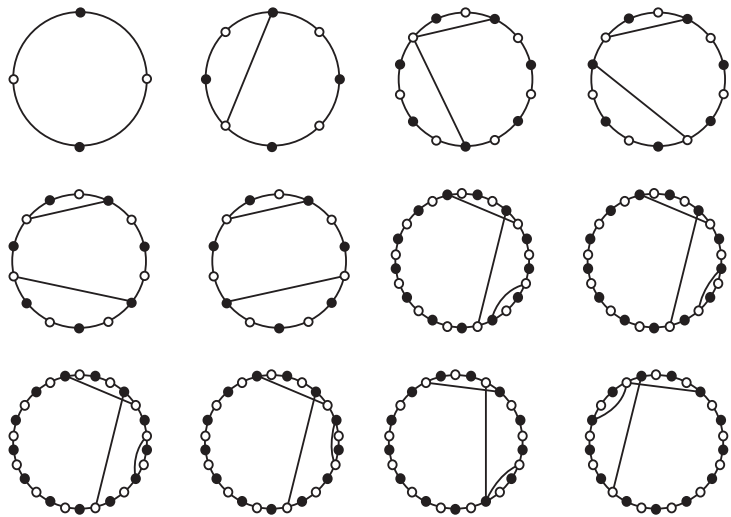}
\caption{The (1)-bipancyclic graphs of order less than 32 }
\label{UPC.le.32}
\end{center}
\end{figure}

\begin{figure}[ht]
\begin{center}
\includegraphics[scale=.600]{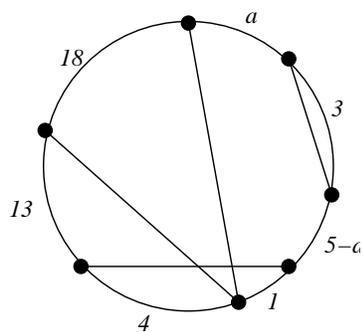}
\caption{The six non-isomorphic (1)-bipancyclic graphs of order 44}
\label{UPC44}
\end{center}\vspace{-.2in}
\end{figure}

Wallis \cite{Wallis} finds all uniquely bipancyclic graphs on at most 30 vertices.
These graphs are displayed in Figure \ref{UPC.le.32}.
In \cite{KPWW}, using computer programs, khodkar et al. show that (see Figure \ref{UPC44}):
\begin{enumerate}
\item If $32\leq v\leq 56$, and $v\neq 44$, then there are no UBPC graphs of order $v$;
\item There are precisely six non-isomorphic UBPC graphs of order $44$ (see Figure \ref{UPC44});
\item There are no other UBPC graphs with $v$ vertices and $v+m$ edges for $0\leq m \leq 5$.
\end{enumerate}

\begin{figure}[ht]
\begin{center}
\begin{tabular}{ccc}
\includegraphics[scale=.500]{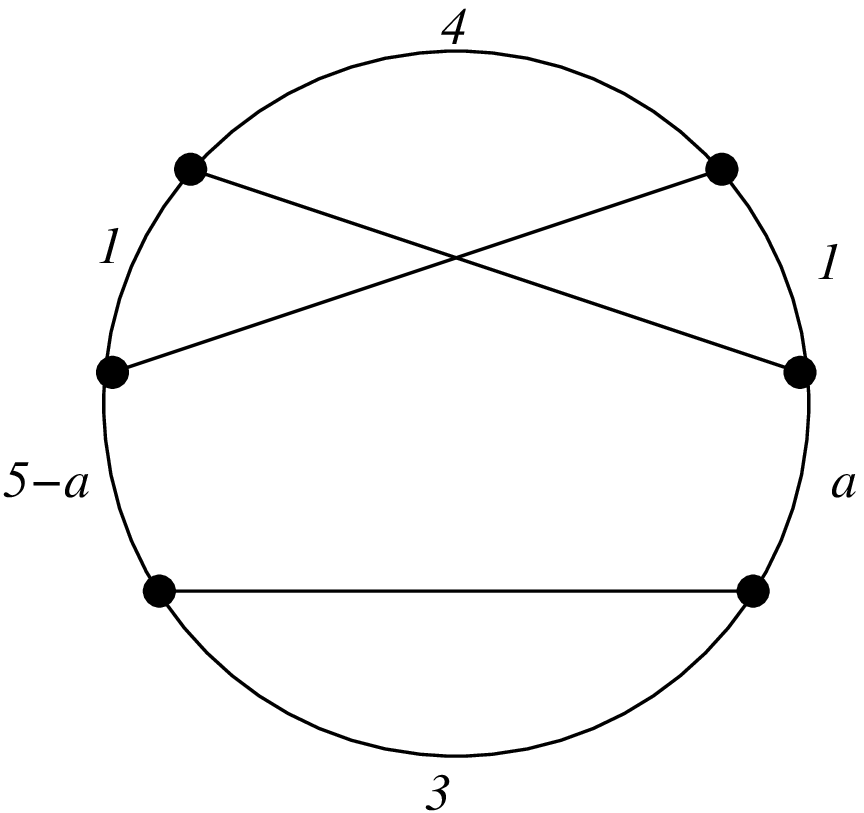}&&
\includegraphics[scale=.500]{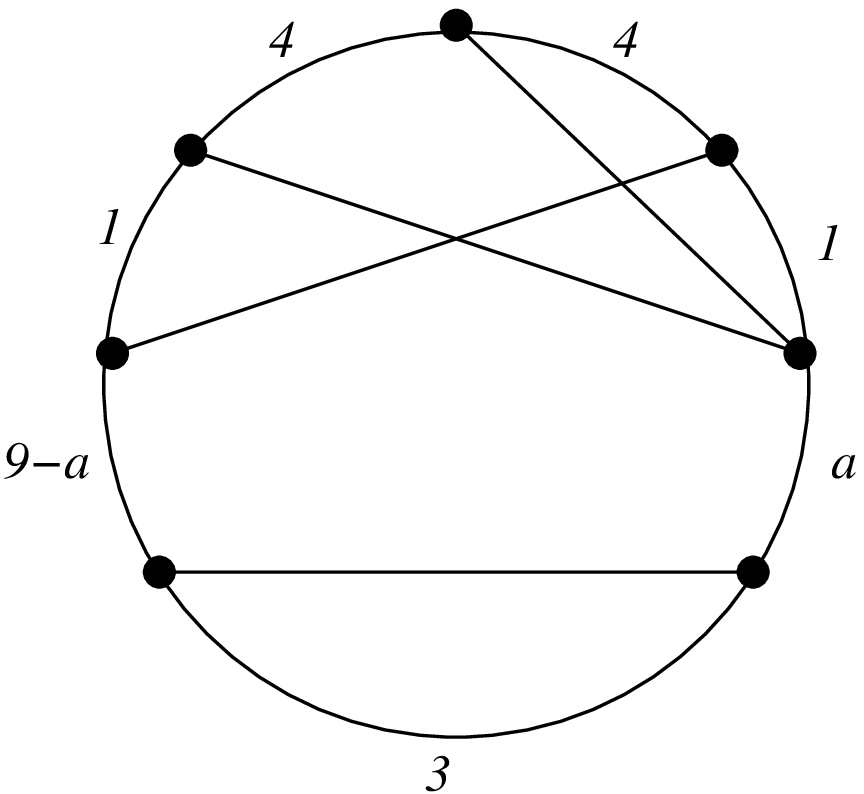}\\
order=14&& order=22\\
3 non-isomorphic && 10 non-isomorphic\\
\end{tabular}
\caption{The (2)-bipancyclic graphs with at most four chords}
\label{2bipan}
\end{center}
\end{figure}

Using computer search we obtain:

\begin{theorem}\label{(r)-bipancyclic}
\begin{enumerate}
\item The graphs displayed in Figure \ref{2bipan} are the only (2)-bipancyclic graphs
with $v$ vertices and $v+m$ edges, where $0\leq m\leq 5$.

\item There is no $(r)$-bipancyclic graph, $r\geq 3$, with $v$ vertices and $v+m$ edges,
where $0\leq m\leq 5$.
\end{enumerate}
\end{theorem}

\section{Oddly $(r)$-bipancyclic graphs}

In this section we initiate the study of oddly $(r)$-bipancyclic graphs. These are the bipartite graphs
with odd number of vertices $v$ and minimum degree at least 2 which have precisely $r$ cycles of
every even length from 4 to $v-1$.

\begin{figure}[ht]
\begin{center}
\begin{tabular}{ccc}
\includegraphics[scale=.380]{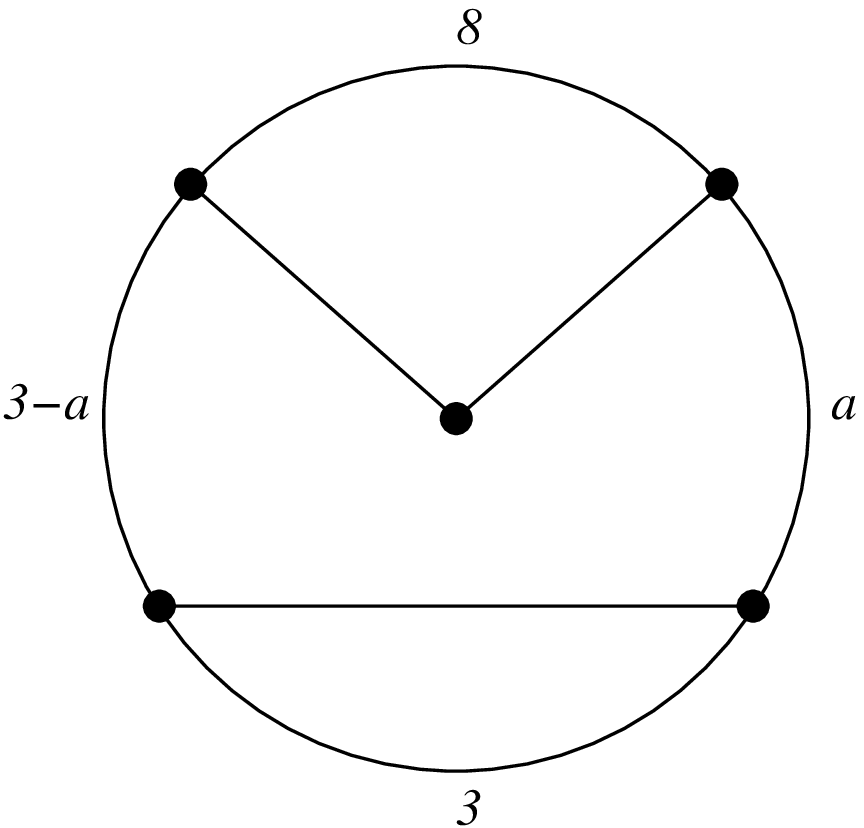}&
\includegraphics[scale=.380]{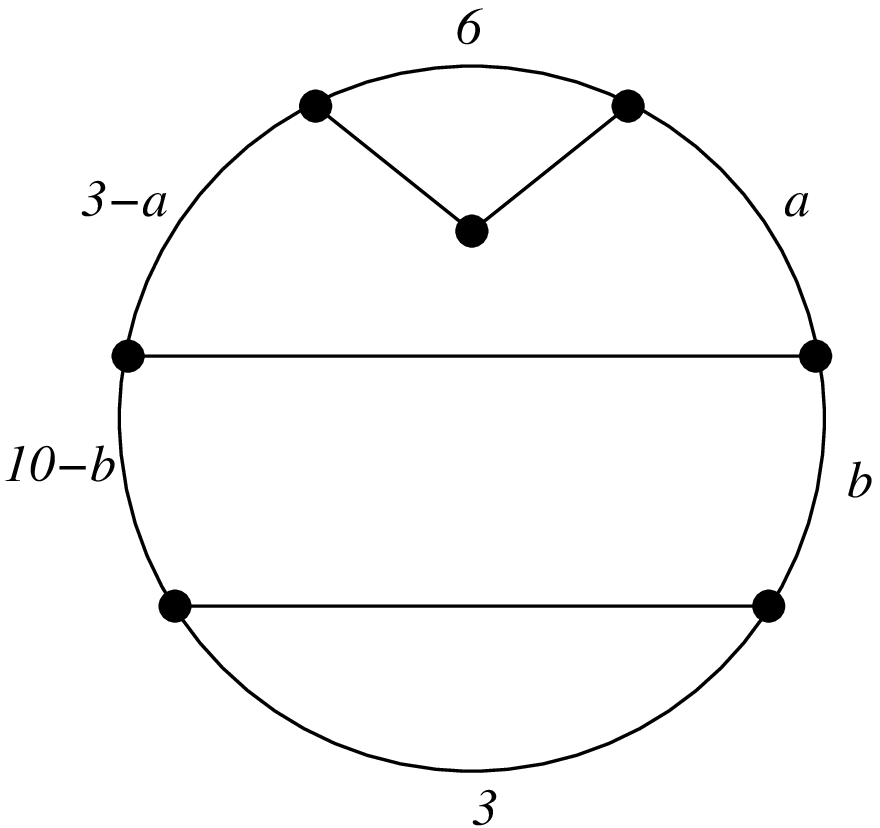}&
\includegraphics[scale=.380]{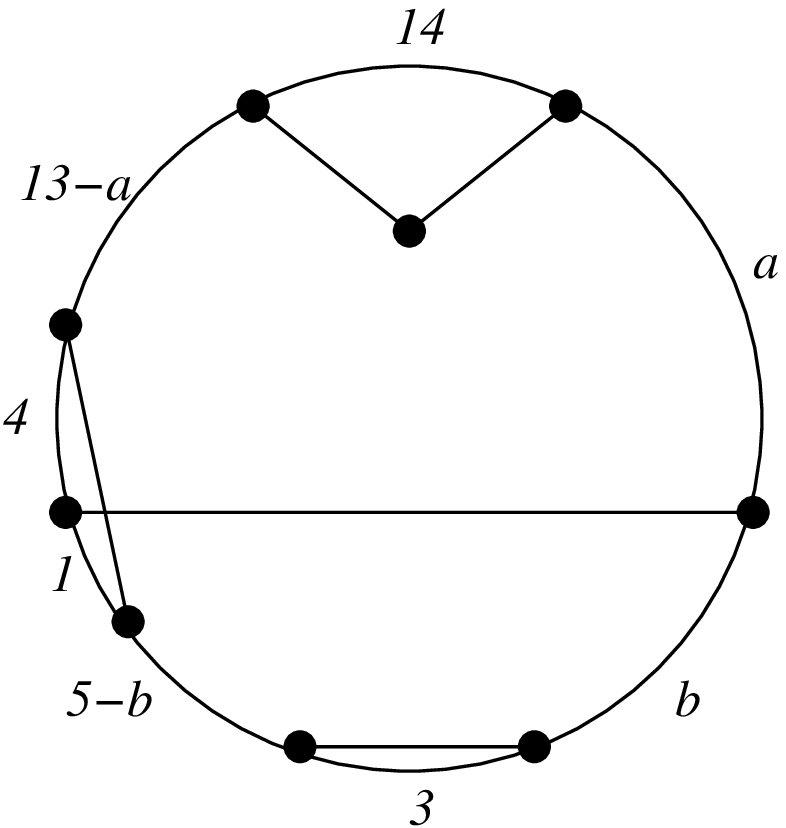}\\
order=15 & order=23 & order=41\\
2 non-isomorphic & 22 non-isomorphic & 84 non-isomorphic\\
\end{tabular}
\caption{The oddly (1)-bipancyclic graphs with $v$ vertices and $(v-1)+m$ edges, where $m\leq 5$.}
\label{uobp123}
\end{center}
\end{figure}

\begin{figure}[ht]
\begin{center}
\begin{tabular}{cc}
\includegraphics[scale=.450]{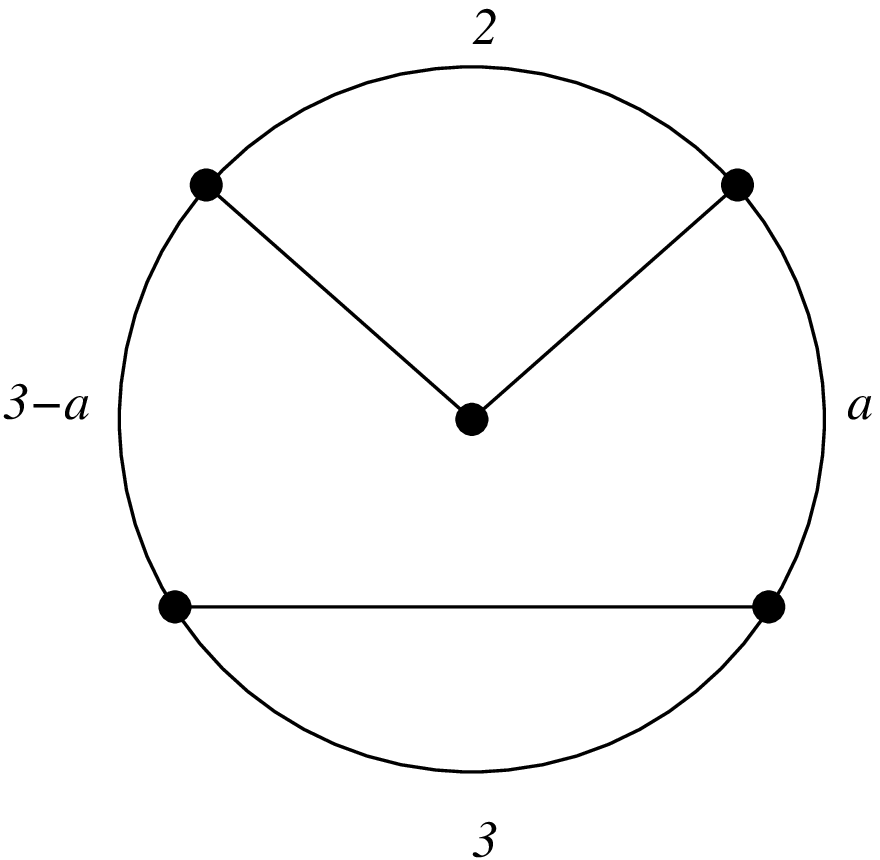}&
\includegraphics[scale=.450]{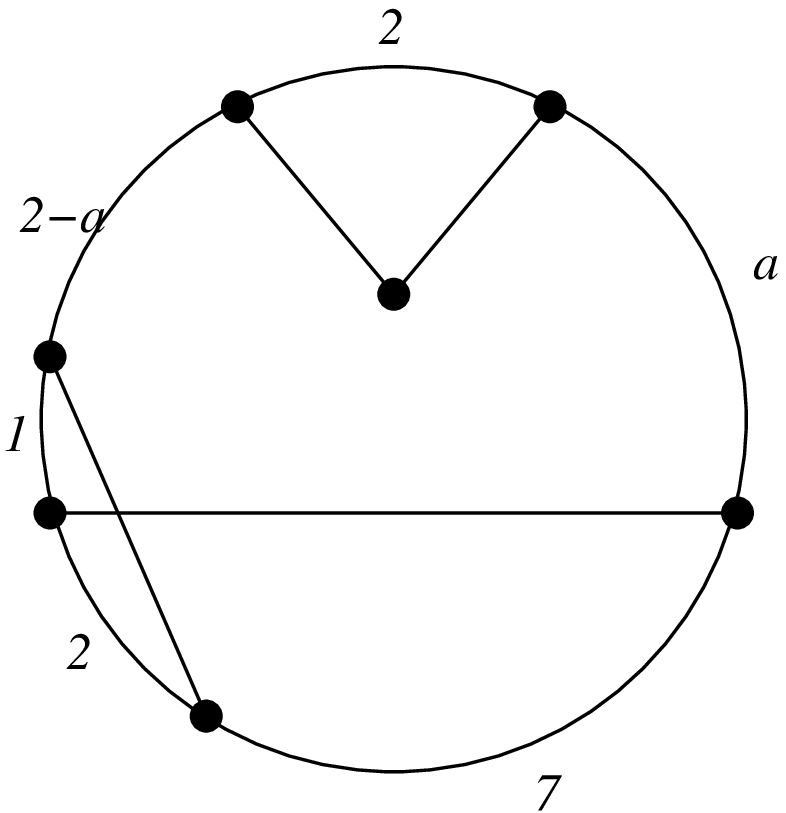}\\
order=9&order=15\\
2 non-isomorphic & 3 non-isomorphic\\
\includegraphics[scale=.450]{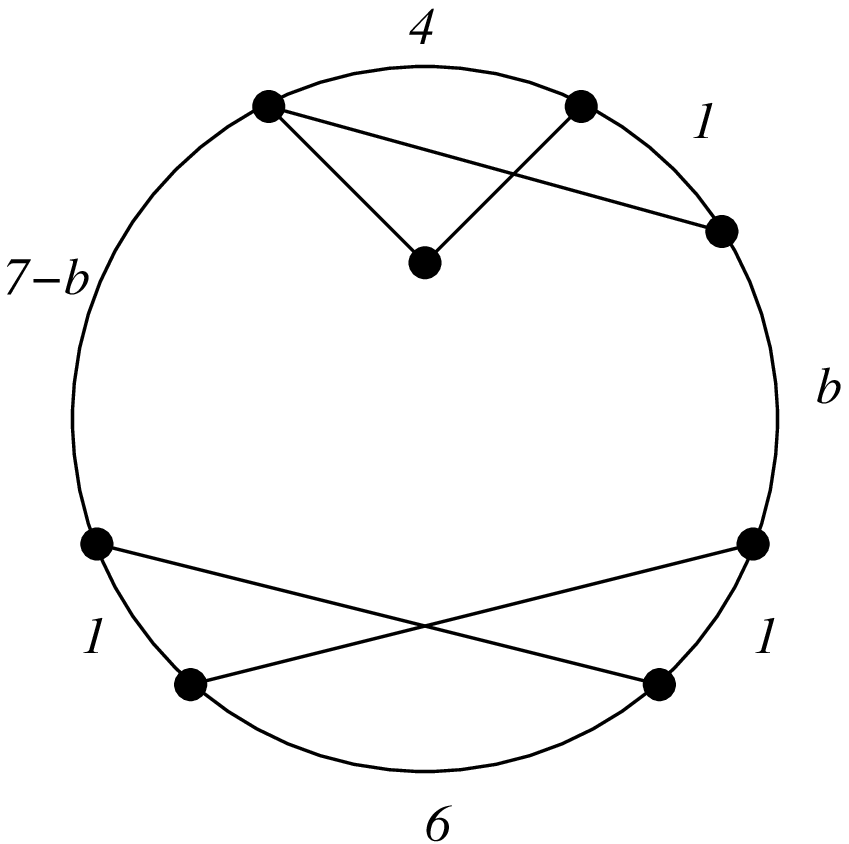}&
\includegraphics[scale=.450]{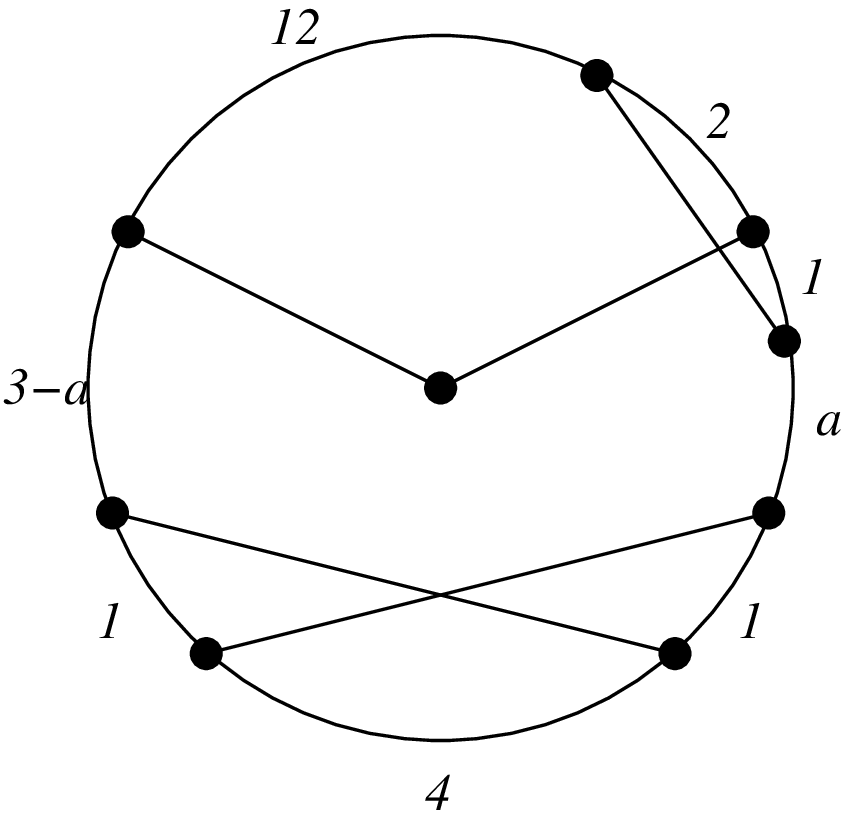}\\
order=21&order=25\\
8 non-isomorphic & 4 non-isomorphic \\

\end{tabular}
\caption{The oddly (2)-bipancyclic graphs with $v$ vertices and $(v-1)+m$ edges, where $m\leq 5$.}
\label{2obp}
\end{center}
\end{figure}

\begin{figure}[ht]
\begin{center}
\begin{tabular}{ccc}
\includegraphics[scale=.400]{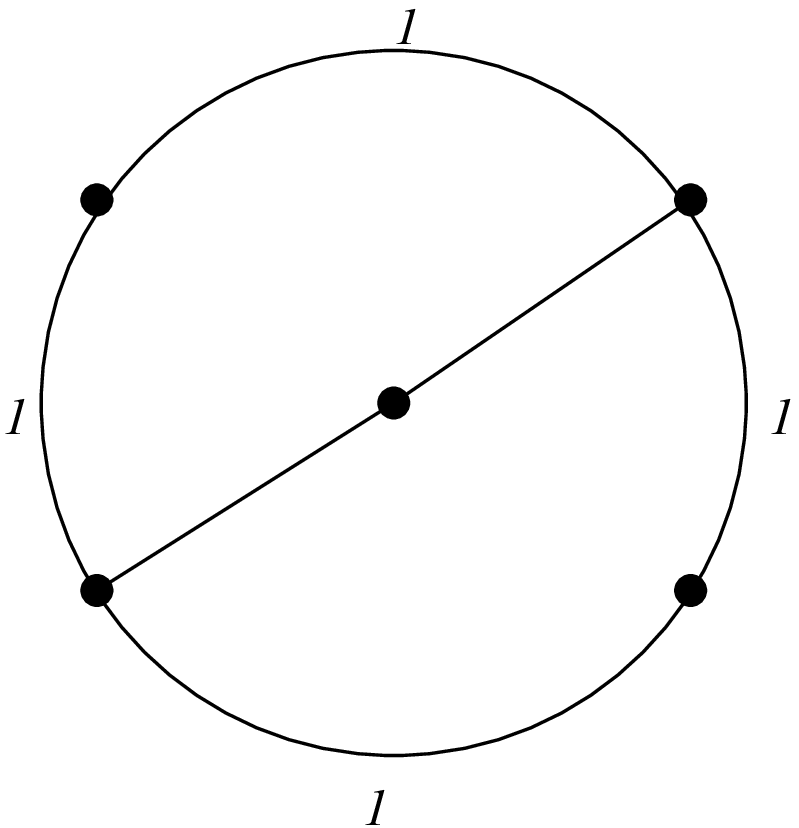}&
\includegraphics[scale=.400]{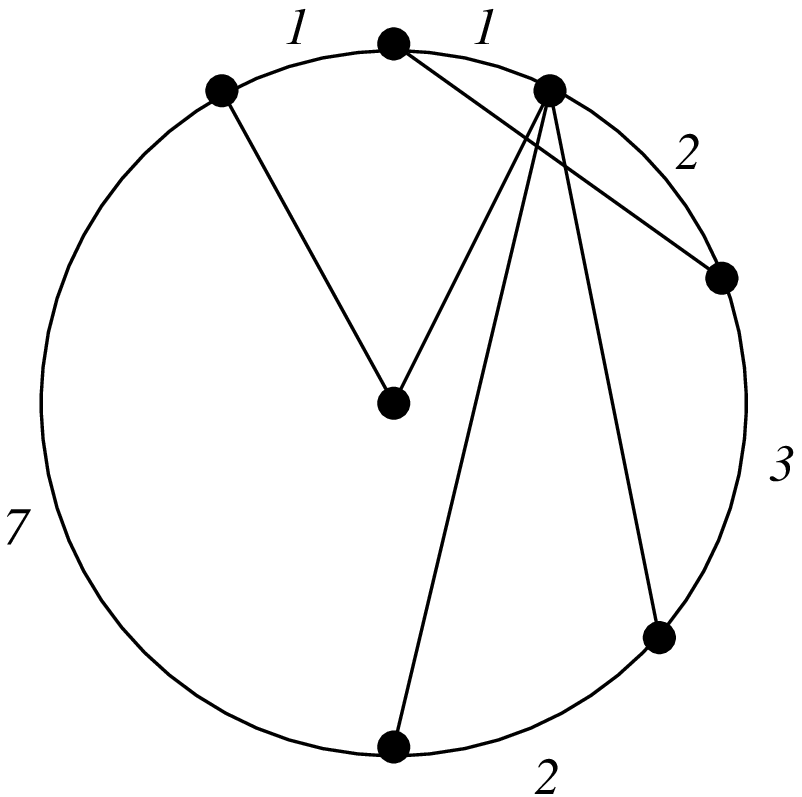}&
\includegraphics[scale=.400]{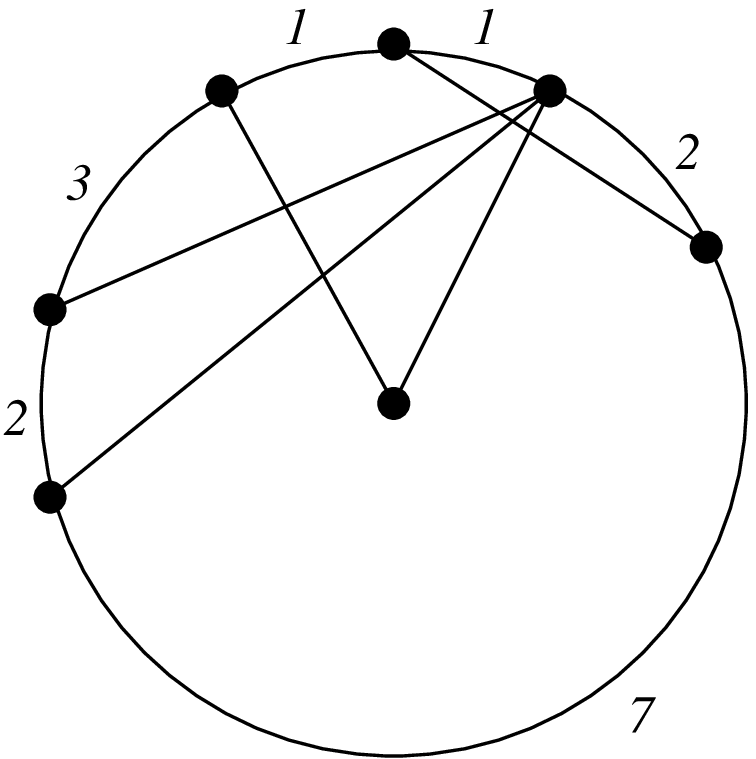}\\
order=5&order=17&order=17\\
\end{tabular}
\caption{The oddly (3)-bipancyclic graphs with $v$ vertices and $(v-1)+m$ edges, where $m\leq 5$.}
\label{3obp}
\end{center}
\end{figure}

\begin{figure}[ht]
\begin{center}
\begin{tabular}{ccc}
\includegraphics[scale=.400]{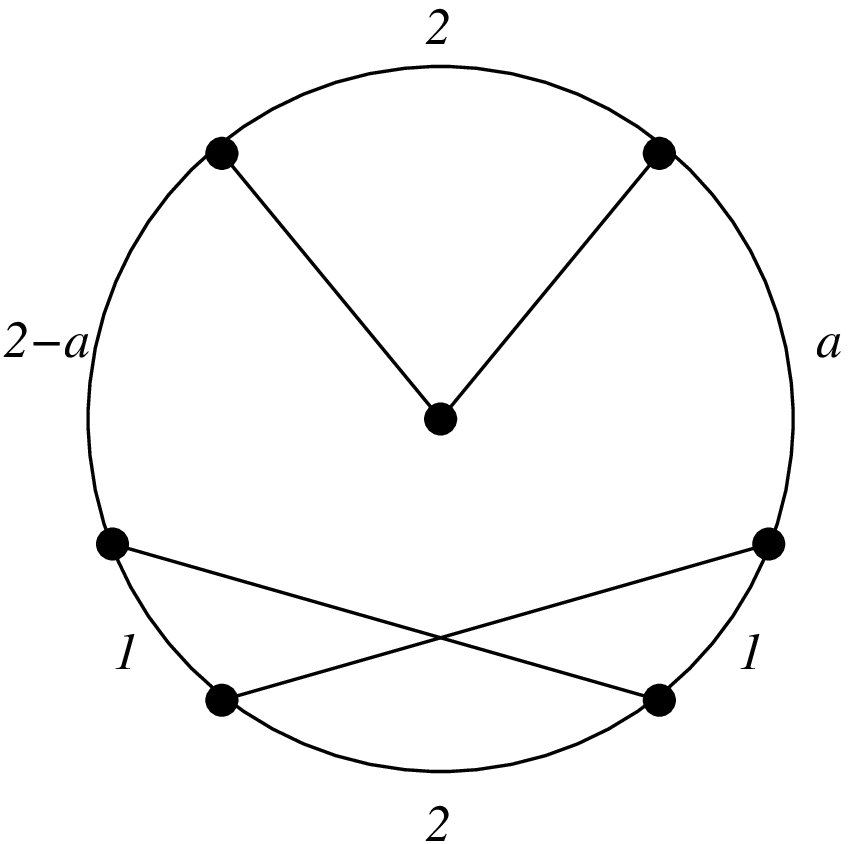}&&\includegraphics[scale=.400]{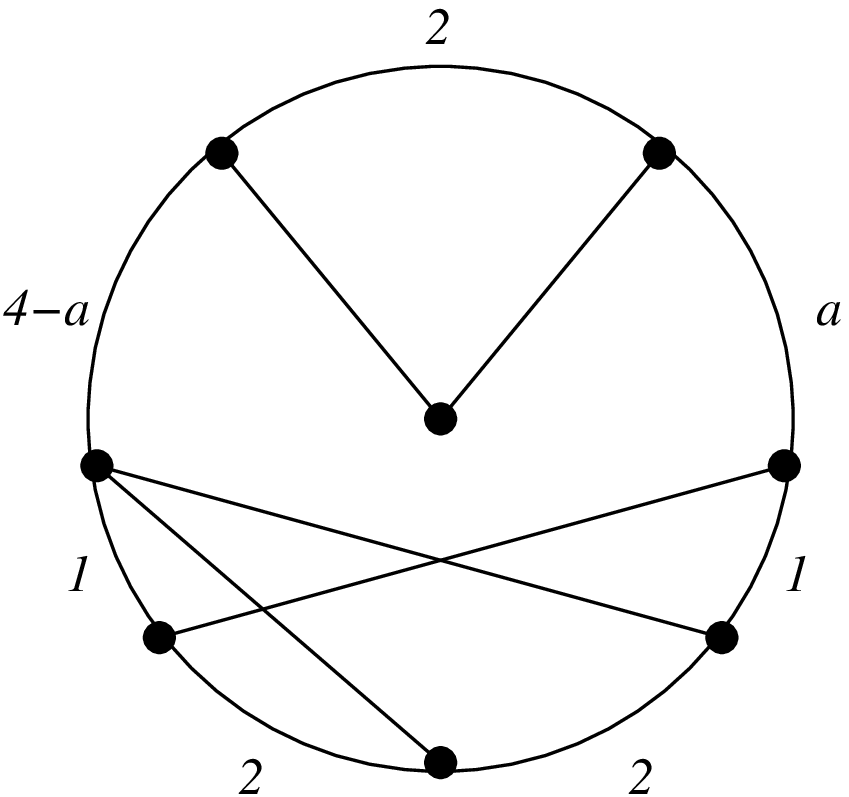}\\
order=9&&order=13\\
2 non-isomorphic && 5 non-isomorphic\\
\end{tabular}
\caption{The oddly (4)-bipancyclic graphs with $v$ vertices and $(v-1)+m$ edges, where $m\leq 5$.}
\label{4obp}
\end{center}
\end{figure}

Our computer search shows that:

\begin{theorem}\label{oddly (r)-bipancyclic}
The graphs displayed in Figures \ref{uobp123}, \ref{2obp}, \ref{3obp} and \ref{4obp}
are the only oddly $(r)$-bipancyclic graphs
with $v$ vertices and $(v-1)+m$ edges, where $2\leq m\leq 5$ and $r\geq 1$.
\end{theorem}

\noindent {\bf Acknowledgement:} The authors would like to thank Alexander Clifton for his
suggestion to study oddly $(r)$-bipancyclic graphs.


\end{document}